\documentclass[journal]{IEEEtran}
\usepackage{amsfonts}
\usepackage{bm}
\usepackage{amsmath,amsfonts,amsthm,amssymb}
\usepackage{setspace}
\usepackage{Tabbing}
\usepackage{fancyhdr}
\usepackage{lastpage}
\usepackage{extramarks}
\usepackage{chngpage}
\usepackage{algorithmic}
\usepackage{soul,color}
\usepackage{epsfig}
\usepackage{graphicx,float,wrapfig,subfigure}
\usepackage{dsfont}
\usepackage{longtable}
\usepackage{wrapfig}
\usepackage{stfloats}
\usepackage{cite}
\usepackage{array}
\usepackage{multirow}{\tiny }
\usepackage{multicol}
\usepackage{colortbl}
\usepackage{tabularx}
\usepackage{mdwmath}
\usepackage{mdwtab}
\usepackage{color}
\usepackage{verbatim}
\usepackage{amsmath}
\usepackage{tikz}
\usetikzlibrary{calc}
\usepackage{flowchart}
\usetikzlibrary{shapes.geometric, arrows}
\usepackage{overpic}
\usepackage{epstopdf}
\usepackage{url}
\usepackage[colorlinks,linkcolor=black,anchorcolor=black,citecolor=black,urlcolor=black]{hyperref} 
\usepackage[flushleft]{threeparttable}
\usepackage{enumitem}

\allowdisplaybreaks
\newcounter{remark}
\newenvironment{remark}[1][]{\refstepcounter{remark}\par
	\textit{Remark~\theremark #1:} \rmfamily}

\newcounter{assumption}
\newenvironment{assumption}[1][]{\refstepcounter{assumption}\par
	\textbf{[A\theassumption #1]:} \rmfamily}

\newcounter{proposition}
\begin{document}

\title{Joint Fleet Sizing and Charging System Planning for Autonomous Electric Vehicles}

\author{
Hongcai~Zhang,~\IEEEmembership{Member,~IEEE,}
Colin J. R. Sheppard,
Timothy E. Lipman,
and~Scott~J.~Moura,~\IEEEmembership{Member,~IEEE}

\thanks{
	This work was funded by the U.S. Department of Energy Vehicle Technologies Office under Lawrence Berkeley National Laboratory Contract No. DE-AC02-05CH11231.
	
	H. Zhang, and S. J. Moura  are with the Department of Civil and Environmental Engineering, University of California, Berkeley, Berkeley, CA 94720 USA. S. J. Moura is also with the Smart Grid and Renewable Energy Laboratory, Tsinghua-Berkeley Shenzhen Institute, Shenzhen, 518055, P.~R.~China (email:  zhang-hc13, smoura@berkeley.edu).
	
	C. J. R. Sheppard is with the Energy Technologies Area, Lawrence Berkeley National Laboratory, Berkeley CA 94720 USA, and also with the Department of Transportation Engineering, University of California, Berkeley, Berkeley, CA 94720 USA (e-mail: colin.sheppard@lbl.gov).
	
	T. E. Lipman is with the Transportation Sustainability Research Center, University of California, Berkeley, Berkeley, CA 94720 USA (e-mail: telipman@berkeley.edu).
}
}

\maketitle

%\large

\begin{abstract} 
	This paper studies the joint fleet sizing and charging system planning problem for a company operating a fleet of autonomous electric vehicles (AEVs) for passenger and goods transportation. Most of the relevant published papers focus on intracity scenarios and adopt heuristic approaches, e.g., agent based simulation, which do not guarantee optimality. In contrast, we propose a mixed integer linear programming model for intercity scenarios.
	This model incorporates comprehensive considerations of 1) limited AEV driving range; 2) optimal AEV routing and relocating operations; 3) time-varying origin-destination transport demands; and 4) differentiated operation cost structure of passenger and goods transportation. The proposed model can be computational expensive when the scale of the transportation network is large. We then exploit the structure of this program to expedite its solution. 
	Numerical experiments are conducted to validate the proposed method. Our experimental results show that AEVs in passenger and goods transportation have remarkable planning and operation differences. We also demonstrate that intelligent routing and relocating operations, charging system and vehicle parameters, e.g., charging power, battery capacity, driving speed etc., can significantly affect the economic efficiency and the planning results of an AEV fleet.
\end{abstract}


\section{Introduction}
\IEEEPARstart{A}{utonomous}  driving is believed to be a disruptive technology that will transform our transportation system in the near future. Autonomous vehicles (AVs) that transport passengers or goods without human intervention will not only free human drivers from burdensome driving labor, but also promote transportation accessibility \cite{Meyer2017}, cut down mobility costs \cite{Bosch2017,HYLAND2018278,FS_Moreno2018}, enhance energy efficiency, and reduce greenhouse gas emission \cite{Greenblatt2015,mahmassani201650th,Milakis2017,Yi2018}. When AVs are electrified, which we refer to as autonomous electric vehicles (AEVs), then the last two aforementioned advantages will be further enhanced, particularly if the electricity is supplied from clean energy (e.g., renewable power generation). 

Passenger (e.g., ride-hailing) and goods transportation are clear initial markets for AEV fleets. Specifically, they feature high utilization levels and planned routes, which can exploit the aforementioned advantages of AEVs\cite{Lam2015,AEVSharing_Stocker2017}. Hence, it is soon possible for transportation network companies, e.g., Uber and Lyft, or logistics companies, e.g., UPS and DHL, to operate a fleet of AEVs in their businesses. 

In this paper, we focus on the planning problem of AEVs for these businesses in intercity transportation. We envision that, with AEVs, transportation network companies may also launch services to satisfy passenger transportation demands across cities, e.g., autonomous electric intercity buses. As a result, they may face similar problems with logistics companies. We will uniformly refer to them as AEV operating companies. 

An AEV operating company needs to solve the following problems before launching its business:
\begin{enumerate}
	\item Fleet sizing, i.e., size an AEV fleet to satisfy given transport demands. On one hand, though vehicle automation allows companies to hire less human drivers to save labor costs, AVs could be very expensive especially in the early stage of commercialization. On the other hand,  high automation allows the company to dispatch and route AEVs more efficiently leading to much higher vehicle utilization. Thus, a company needs to optimally design its AEV fleet size to fully exploit AEVs' potential and reduce its costs considering future strategic operations.
	\item Charging system planning, i.e., site and size EV chargers to satisfy charging demands. Compared with traditional internal combustion vehicles, EVs generally have much higher fuel efficiency, thus will significantly reduce a company's fuel costs. However, the company may also need to invest in sufficient charging infrastructure to charge the AEVs\cite{neubauer2014impact}. This can be vitally important for the early stage of transportation electrification when public fast chargers are not quite popularized, especially on intercity corridors. 
\end{enumerate}
The above two problems are highly coupled together because the planning of charging systems can impact the utilization of AEVs: 1) charging can lead to significant ``down time'';\footnote{Because the rated charging power is not high enough.} 2) AEVs may detour to get charged when chargers are not available on their preferred paths or the charging prices elsewhere are cheaper so that vehicle miles traveled may increase. 

As a result, it is necessary for a company to jointly design its fleet size and charging system according to the expected transportation demands. It is also important to carefully consider the routing and charging operations in the design solution.

Both AV fleet sizing and EV charging system planning have been active areas of research for years. Many researchers have studied these two problems separately, summarized as follows.

High automation of AVs means they have the potential to satisfy more demands than the same number of human-driven vehicles. Hence, many researchers have studied the AV fleet sizing problem to evaluate their economic advantages. 
Boesch et. al. \cite{boesch2016autonomous} proposed an agent-based simulation approach for the fleet sizing problem of shared-use AVs based on which the authors evaluated how the AVs' service level could affect the required fleet size. 
Alonso-Mora et. al. \cite{Alonsomora2017} developed a dynamic
trip-vehicle assignment strategy for autonomous ride-sharing services and studied the trade-off between the fleet size and the performance of the ride-sharing service. 
Fagnant et. al. \cite{FS_Fagnant2018} developed an agent- and network- based simulation framework that considers dynamic ride-sharing and vehicle relocation to estimate fleet size and operator profitability for given demands. 
Vazifeh et. al. \cite{FS_Nature_Vazifeh2018} provided a network-based solution to the fleet sizing problem to determine the minimum number of vehicles needed to serve given trips without incurring any delay to the passengers.

Limited driving range is the major hurdle for EV application. The charging system planning problem has been extensively studied in the literature. Generally, the published approaches can be divided into three categories: 1) computational geometry based approaches that assume charging demands occur on geographical nodes and ignore transportation network structure \cite{Plan_T_Xu2013,Plan_T_Cavadas2015,Plan_Zhang2015,Faridimehr2018}; 2) Origin-Destination (OD) flow based methods that explicitly describe transportation network constraints, i.e., driving range constraints \cite{FRLM_Kuby2005,FRLM_MirHassani2013,Plan_T_Mak2013,Plan_T_You2014,Plan_T_Riemann2015,FRLM_Chung2015,XWang2017,YXiong2018}; 3) simulation- or data- driven approaches \cite{Plan_T_Cai2014,Plan_T_Shahraki2015}, which adopt agent-based simulation or real-world data to estimate future EV charging behaviors. Because that EV charging systems are coupling power and transportation networks together, many researchers have proposed multidisciplinary approaches to plan charging systems in coupled networks \cite{Plan_TE_He2013,Plan_TE_Wang2013,Plan_TE_Sadeghi2014, Plan_TE_Yao2014,Plan_TE_Luo2015,Plan_TE_Xiang2016}. 

However, there are only few papers that have addressed the joint AEV fleet sizing and charging system planning problem.
Hiermann et. al. \cite{FSCS_Hiermann2016} studied the EV fleet sizing and routing problem including the choice of recharging times and locations. They took the design of EV charging stations as given. 
Chen et. al. \cite{FS_Chen2016} proposed an agent-based simulation platform to size shared-use AEVs. Based on the platform, they evaluated various charging infrastructure investment decisions. 
Bauer et. al.\cite{Bauer2018} deigned an agent-based simulation framework to estimate required fleet size and most economic charging systems to satisfy taxi-hailing demands. Dandl et. al. \cite{Dandl2018} compared fleet size requirement and profitability of an AEV taxi fleet with an existing free-floating carsharing system based on a data-driven simulation approach.

None of the above papers have addressed the intercity scenarios of AEV application. Most of the proposed approaches are based on simulation and heuristic optimization techniques that do not guarantee optimality. Besides, most of the papers focus on passenger transportation while AEVs may behave quite differently when they are used for goods transportation: For the former, it is necessary to drop the passenger to the destination in time because the passenger's time is quite valuable; However, for the latter, an AEV has lower time cost so that it will be willing to detour more so that it can get cheaper electricity or avoid congested charging stations\cite{Yu2018}. A comparison on the operation costs between passenger and goods transportation is given in Table~\ref{tab_application}. The \textit{loaded} trips are with passengers or goods, while the \textit{relocating} ones are empty AEVs relocating from one location to another.

To address the aforementioned problem, we propose a joint fleet sizing and charging system planning model to help AEV operating companies to right-size their AEV fleets and design charging infrastructure at least cost while meeting their transport demands in intercity corridors. This paper advances the relevant literature by the following contributions:
\begin{enumerate}
	\item A joint AEV fleet sizing and charging system planning model that simultaneously optimizes the fleet size and charging system to satisfy a mobility on demand system is developed. The planning model takes the strategic routing and relocating of AEVs into consideration so that it can balance the investment costs at the planning stage and the operation costs in the future. Furthermore, both passenger and goods transportation can be modeled.
	\item The limited driving range constraints of AEVs are explicitly described by an expanded transportation network model based on OD transport demands. The impact of vehicle charging behaviors on fleet operation and charging system planning can be effectively evaluated. 
	\item The proposed model is a mixed integer linear program that can be computationally expensive when the scale of the transportation network is large. We exploit the structure of the program and propose an efficient algorithm to expedite its solution with guaranteed solution quality.
\end{enumerate}
In addition, numerical experiments are conducted to validate the proposed method. The results show that AEVs in different application scenarios, i.e., passenger and goods transportation, have fundamental planning and operation differences. We also demonstrate that intelligent relocating operation and system and vehicle parameters, e.g., charger power, battery capacity, driving speed etc., all can significantly affect the economic efficiency of an AEV fleet and the planning results.

We briefly introduce the expanded transportation network model in Section \ref{sec_transportation}. Section \ref{sec_model} formulates the complete planning model. We exploit the structure of the model and propose a solution approach in Section \ref{sec_algorithm}. Numerical experiments are presented in Section \ref{sec_case}. Section \ref{sec_conclusion} concludes the paper.

\begin{table}
	\renewcommand{\arraystretch}{1.25}
%	\vspace{-8mm}
	\centering
	\begin{footnotesize}
		\caption{Comparison between passenger and goods transportation}
		\vspace{-2mm}
		\begin{tabular}{ccccc}
			\hline
			{Case}& Trip type&Passenger time& Fuel costs & Maintenance\\
			\hline
			\multirow{2}{*}{Passenger} & Loaded & Yes & Yes & Yes\\
			& Relocating & No & Yes & Yes\\
			\multirow{2}{*}{Goods} & Loaded & No &Yes&Yes\\
			& Relocating & No &Yes&Yes\\
			\hline 
		\end{tabular}\label{tab_application}
	\end{footnotesize}
	\vspace{-4mm}
\end{table}

\section{Expanded Transportation Network Model}\label{sec_transportation}

This section introduces an expanded transportation network model that describes AEVs' driving range constraints. 
This model originates from the flow-refueling location model developed by Kuby and Lim \cite{FRLM_Kuby2005} that assumes AEVs' driving paths are fixed. This proposed model extends it by allowing AEVs to choose different paths to drive from an origin to its destination.
To construct a computational efficient transportation network model with driving range constraints, we assume that:
\begin{assumption}\label{as_charge}
	AEVs will get fully charged at charging stations.\footnote{There is usually significant fixed time cost for an AEV to use a charging station, it will not tend to get partially charged each time.}
\end{assumption} 
\begin{assumption}\label{as_drive}
	The routing of AEVs will not affect traffic conditions, so the driving speed on each arc is exogenously defined.\footnote{This assumption may be realistic when the majority of vehicles on road are highly automatic and intelligent so that the network is rarely congested, or the penetration of AEVs in the vehicles on road is small so that they have negligible impacts on traffic conditions.} %and their routing will not significantly affect traffic congestion and traveling speed
\end{assumption} 

We use a directed graph $G({\mathcal{I}},\mathcal{A})$ to model the transportation network, where $\mathcal{I}$ denotes node set and $\mathcal{A}$ denotes arc set. An \textit{arc} $(i,j)\in\mathcal{A}$ is the road link between two adjacent nodes, $i\in\mathcal{I}$ and $j\in\mathcal{I}$. 
Symbol $l_{ij}$ is the length of arc $(i,j)$.
We use OD pairs to model transport demands. An OD pair (indexed by $g$) is composed by an \textit{origin node} $o_g\in\mathcal{I}$ and a \textit{destination node} $d_g\in\mathcal{I}$. The transport demands in one OD pair $g$ is denoted by tuple $(o_g,d_g,\lambda_{g}^\text{od})$, $\forall g\in\mathcal{G}$. Symbol $\lambda_{g}^\text{od}$ represents the Poisson volume of demands from origin node $o_g$ to destination node $d_g$. 
We call the route that an AEV can choose from an origin to a destination as a \textit{path}.

\subsection{Transportation Network Expansion}
\begin{figure}
%	\vspace{-10mm}
	\centering
	\subfigure[The original network $G({\mathcal{I}},{\mathcal{A}})$]{
		\resizebox{0.425\textwidth}{!}{
			\centering
			\tikzset{>=latex, every picture/.style={line width=0.75pt}}
			\begin{tikzpicture}
			\begin{footnotesize}
			\node[draw,circle,minimum size=4.mm,xshift=-10mm,outer sep=0](s) {o};
			\node[draw,circle,minimum size=4.mm,xshift=00mm,outer sep=0](A)  {1};
			\node[draw,circle,minimum size=4.mm,xshift=15mm,outer sep=0](B) {2};
			\node[draw,circle,minimum size=4.mm,xshift=15mm,yshift=7.5mm,outer sep=0](E) {5};
			\node[draw,circle,minimum size=4.mm,xshift=40mm,outer sep=0](C) {3};
			\node[draw,circle,minimum size=4.mm,xshift=55mm,outer sep=0](D) {4};
			\node[draw,circle,minimum size=4.mm,xshift=65mm,outer sep=0](t) {d};
			\draw[->] (A.east) -- (B.west) node[below=1.5pt, xshift=-5mm] {30 km};
			\draw[->] (B.east) -- (C.west) node[below=1.5pt, xshift=-10mm] {60 km};
			\draw[->] (C.east) -- (D.west) node[below=1.5pt, xshift=-5mm] {30 km};
			\draw[->] (s.east) -- (A.west) node[below=1.5pt, xshift=-2.5mm] {50km};
			\draw[->] (D.east) -- (t.west) node[below=1.5pt, xshift=-2.5mm] {50km};
			
			\draw[->] (A.east) -- (E.west) node[left=1.5pt, yshift=-1mm, xshift=-3mm] {35 km};
			\draw[->] (E.east) -- (C.west) node[right=1.5pt, yshift=7mm, xshift=-14mm] {65 km};
			\end{footnotesize}
			\end{tikzpicture}\label{fig_path1}
	}}
	\subfigure[The expanded network $G(\mathcal{I},\hat{\mathcal{A}})$ (driving range 100 km)]{
		\resizebox{0.45\textwidth}{!}{
			\centering
			\tikzset{>=latex, every picture/.style={line width=0.75pt}}
			\begin{tikzpicture}
			\begin{footnotesize}
			\node[draw,circle,minimum size=4.mm,xshift=-10mm,outer sep=0](s) {o};
			\node[draw,circle,minimum size=4mm,xshift=00mm,outer sep=0](A)  {1};
			\node[draw,circle,minimum size=4.mm,xshift=15mm,outer sep=0](B) {2};
			\node[draw,circle,minimum size=4.mm,xshift=15mm,yshift=7.5mm,outer sep=0](E) {5};
			\node[draw,circle,minimum size=4.mm,xshift=40mm,outer sep=0](C) {3};
			\node[draw,circle,minimum size=4.mm,xshift=55mm,outer sep=0](D) {4};
			\node[draw,circle,minimum size=4.mm,xshift=65mm,outer sep=0](t) {d};
			\draw[->] (s.east) -- (A.west) node[above=1.5pt, xshift=-2.5mm] {};
			\draw[->] (A.east) -- (B.west) node[above=1.5pt, xshift=-5mm] {};
			\draw[->] (B.east) -- (C.west) node[above=1.5pt, xshift=-13mm] {};
			\draw[->] (C.east) -- (D.west) node[above=1.5pt, xshift=-3mm] {};
			\draw[->] (D.east) -- (t.west) node[above=1.5pt, xshift=-2.5mm] {};
			
			\draw[->] (A.east) -- (E.west) node[above=1.5pt, xshift=-3mm] {};
			\draw[->] (E.east) -- (C.west) node[above=1.5pt, xshift=-3mm] {};
			
			\path[every node/.style={font=\sffamily\small}]
			(s) edge[->,bend right,dashed] node [left] {} (B);
			\path[every node/.style={font=\sffamily\small}]
			(C) edge[->,bend left,dashed] node [left] {} (t);
			\path[every node/.style={font=\sffamily\small}]
			(B) edge[->,bend right,dashed] node [left] {} (D);
			\path[every node/.style={font=\sffamily\small}]
			(A) edge[->,bend right,dashed] node [left] {} (C);
			\path[every node/.style={font=\sffamily\small}]
			(s) edge[->,bend left,dashed] node [left] {} (E);
			
			\end{footnotesize}
			\end{tikzpicture}
			\label{fig_path2}
	}} 
	\vspace{-2mm}
	\caption{An example of the transportation network expansion. The numbers in the circles denote the indexes of the transportation nodes. The arrowed lines represent the directed arcs. The distances of the arcs are marked next to them. This network has 5 nodes, ${\mathcal{I}}=\{{o,1,2,3,4,5,d}\}$, and 7 arcs, $\mathcal{A}=\{{(o,1),(1,2),(1,5),(2,3),(3,4),(4,d),(5,3)}\}$. 
		The original transportation network $G(\mathcal{I},\mathcal{A})$ in  Fig.~\ref{fig_path1} is expanded into a new network  $G(\mathcal{I},\hat{\mathcal{A}})$ in Fig.~\ref{fig_path2} to describe PEV driving range constraints, in which $\hat{\mathcal{A}}=\{(o,1), (o,2),(o,5),(1,2),(1,3), (1,5),(2,3), (2,4), (3,4),$ $(3,d),(4,d),(5,3)\}$.
	}
	\label{fig_network_expansion}
	\vspace{-5mm}
\end{figure}
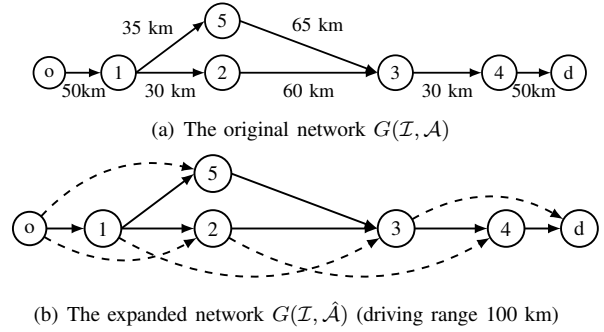

We take a simple transportation network $G({\mathcal{I}},\mathcal{A})$ in Fig.~\ref{fig_path1} as an example. It only has a single OD pair, $g=(o,d)$. There are two paths that the AEVs can travel from origin node $o$ to destination node $d$, i.e., path $(o1534d)$ and path $(o1234d)$. 
The AEV transport demand in this network is $(o,d,\lambda_g^\text{od})$. The AEVs' driving range is $R=100$ km. The average driving speeds on different arcs are the same. 

We can expand the network $G(\mathcal{I},\mathcal{A})$ in  Fig.~\ref{fig_path1} into a new network  $G(\mathcal{I},\hat{\mathcal{A}})$ in Fig.~\ref{fig_path2} by the following steps:
\begin{enumerate}
	\item Connect any two nodes, say $i$ and $j$, in $\mathcal{I}$, by a pseudo arc if node $j$ can be reached from node $i$ after a single charge. Let each new arc's direction be consistent with the traffic flow direction of the original network.
	\item Let the lengths of the pseudo arcs defined in step 1) be equal to those of the corresponding two nodes' shortest paths. For example, $l_{o2}=l_{o1}+l_{12}$. For arc $(1,3)$, an AEV may drive through it by path $(123)$ or $(153)$. However, under Assumption [A2], an AEV driving from node $1$ to node $3$ without getting charged between them will only choose the shortest path $(123)$ to minimize its costs. Therefore, we can let $l_{13}=l_{12}+l_{23}$. 
	\item Let set $\hat{\mathcal{A}}$ be the union of $\mathcal{A}$ and the pseudo arcs added in step 1). Then, we have the expanded network $G(\mathcal{I},\hat{\mathcal{A}})$.
\end{enumerate}

In the expanded network $G({\mathcal{I}},\hat{\mathcal{A}})$, each path from $o$ to $d$ characterizes a feasible solution for an AEV with driving range $R$ to drive from node $o$ to $d$ on the condition that whenever it runs across a node in $\mathcal{I}$ it get charged. For example, an AEV can travel through path \{$o24d$\} if it gets charged at nodes $2$ and $4$; and it can also travel through path \{$o53d$\} if it get charged at nodes $3$ and $5$. In summary, this expanded transportation network model incorporates AEV driving range constraints.

Note that under assumption [A\ref{as_drive}], when an AEV does not get charged between two nodes $i$ and $j$, it will only choose the shortest path (with minimum time or costs) to drive from $i$ to $j$, rendering a unique pseudo arc $(i,j)$. Thus, given the original transportation network, the expanded network can be uniquely determined \textit{a priori} off-line and the cardinality of the arc set in the expanded network is limited by $|I|^2$. 
Furthermore, the energy consumption and driving time on each arc can also be estimated \textit{a priori} based on the arc length. 

\subsection{Traffic Flow Continuity Constraints}
For each transport demand $(o_g,d_g,\lambda_{g}^\text{od})$, $\forall g\in\mathcal{G}$, it may traverse any arc and node in the transportation network. The inflow and outflow should balance at each node. Besides, the traffic flow from each origin node and each destination node should be equal to the total transport demand. These constraints can be represented as follows:
\begin{align}
&\sum_{\{j|(i,j)\in \hat{\mathcal{A}}\}}\lambda_{g,ij}-\sum_{\{k|(k,i)\in \hat{\mathcal{A}}\}}\lambda_{g,ki}=
\left\{
\begin{array}{l}
\lambda_{g}^\text{od},~\text{if}~i=o_g\\
-\lambda_{g}^\text{od},~\text{if}~ i=d_g\\
0,~\text{if}~i\neq o_g,d_g
\end{array}\right.,\notag\\
&\qquad\qquad\qquad\qquad\forall g\in \mathcal{G}, \forall i \in \mathcal{I},\label{eqn_chap2_subchap1_1}\\
&\lambda_{g,ij}\geq 0, \qquad \forall g \in \mathcal{G}, \forall (i,j)\in \hat{\mathcal{A}},\label{eqn_chap2_subchap1_2}
\end{align}
where, $\lambda_{g,ij}$ is the portion of AEV flow driving on arc $(i,j)$ corresponding to demand $(o_g,d_g,\lambda_{g}^\text{od})$. Note that the AEVs in $\lambda_{g,ij}$ also get charged in the charging station at node $j$.

For brevity, we can also represent the above constraints as:
\begin{align}
&\bm{A}{{\bm{\lambda}}}_g={\bm{\lambda}}_g^\text{od},\qquad\forall g\in \mathcal{G},\label{eqn_TS1}\\
&{{\bm{\lambda}}}_g\geq 0,\qquad\forall g\in \mathcal{G},\label{eqn_TS2}
\end{align}
where, $\mathbf{A}$ is the node-arc incidence matrix of $G(\mathcal{I},\hat{\mathcal{A}})$; ${\bm{\lambda}}_g^\text{od}$ is the vector of nodal transport demands injections. The elements of ${\bm{\lambda}}_g^\text{od}$ are $\lambda_{g}^\text{od}$ at the origin node; $-\lambda_{g}^\text{od}$ at the destination node; zero at any other node. ${{\bm{\lambda}}}_g$ is the vector of AEV flows on arcs of the expanded network; its elements are $\lambda_{g,ij}, \forall (i,j)\in\mathcal{A}$.

\section{Fleet Sizing and Charging System Planning}\label{sec_model}
In this section, we formulate a joint AEV fleet sizing and charging system planning model, which optimally determines the AEV fleet size and the capacities of AEV charging stations in the transportation network. The nomenclature of this model is summarized in Table \ref{tab_nomenclature1}. We consider time-varying traffic flow and adopt tuple $(o_g,d_g,\lambda_{g, t}^\text{od})$ to denote the transport demands, which means there are $\lambda_{g, t}^\text{od}$ demands that need to depart from origin $o_g$ during hour $t$ to destination $d_g$.

The proposed model considers ``relocation''  of autonomous AEVs to fully utilize their automation potential:
\begin{enumerate}
	\item When an AEV arrives at its destination, it can be used to satisfy another demand that originates from its destination right after it gets fully charged. 
	\item An idle AEV can also be moved to other transportation nodes to satisfy future demands there, which is referred to as ``relocation" in this paper. 
\end{enumerate}
Hereafter, we call the driving AEVs \textit{loaded} ones if they are with passengers or goods and \textit{relocating} ones otherwise.

To construct a computational efficient model and avoid comprehensive and intractable modeling of multiple-trip AEV driving range constraints, we make the following assumption:

\begin{assumption}\label{as_rebalance}
	An AEV will get fully charged when it arrives at a destination before it can be scheduled for the next trip.\footnote{This assumption can be true when each trip is long enough in intercity corridors. With this assumption, we will make a conservative planning result so that the system can be more robust.}
\end{assumption} 

We begin with a base model for passenger transportation, in which a loaded AEV has significant passenger time costs while a relocating one does not as indicated in Table~\ref{tab_application}. Then, we will show how the model can be modified when the cost factors for loaded and relocating AEVs are similar in goods transportation scenarios. 

\begin{table}
		\renewcommand{\arraystretch}{1.15}
	%		\vspace{-2mm}
	\begin{small}
		\caption{Nomenclature of the planning model }\label{tab_nomenclature1}
		\begin{tabular}{p{0.8cm}p{7cm}}
			\hline
			\multicolumn{2}{c}{\textbf{Indices/sets}}\\
			$i,j/\mathcal{I}$ & Index/set of transportation nodes.\\
			$(i,j)/\mathcal{A}$ & Index/set of transportation arcs. $\hat{\mathcal{A}}$ is the corresponding set of the expanded transportation network.\\
			$g/ \mathcal{G}$ & Index/set of OD pairs. \\
			$o_g/d_g$ & Origin/destination node of OD pair $g$. $o_g/d_g\in\mathcal{I}$.\\
			$q$/$\mathcal{Q}_g$,$\mathcal{Q}_g^\text{r}$ & Index/set of paths. $\mathcal{Q}_g$ and $\mathcal{Q}_g^\text{r}$ are sets of paths in OD pair $g$ for loaded and relocating AEVs, respectively. \\
			$\tau$& Index of time intervals, $\Delta \tau$=1 hour in this paper, $0$ and $T$ are the first and last time intervals, respectively.\\
			\hline
			\multicolumn{2}{c}{\textbf{Matrix}}\\
			$\mathbf{A}$ &  Node-arc incidence matrix of network $G(\mathcal{I},\hat{\mathcal{A}})$ \\
			$\mathbf{B}_g$/$\mathbf{B}_g^\text{r}$ & Fundamental path-arc incidence matrix for loaded/relocating AEVs correlated to OD pair $g$ of network $G(\mathcal{I},\hat{\mathcal{A}})$. \\
			\hline
			\multicolumn{2}{c}{\textbf{Parameters}}\\
			$c_j^\text{e}$ & Per-unit electricity cost at location $j$, in \$/kWh.\\
			$c^\text{ev}$ & Per-unit cost for purchasing one AEV, in \$.\\
			$c^\text{m}$ & Per-unit AEV maintenance cost, in \$/km.\\
			$c_{g,ij}^\text{oper}$ & Total operation cost (for time, electricity and maintenance) when an AEV in OD pair $g$ drives on arc $(i,j)$. \\
			$c^\text{sp}$ & Per-unit cost for installing one charger, in \$.\\
			$c^\text{t}$ & Per-unit passenger time cost, in \$/hour.\\

			$l_{ij}$ & Length of arc $(i,j)$.\\
			$p^\text{sp}$ & Rated charging power of chargers, in kW.\\
			$t_{g,ij}$ & Total time (driving \& charging) that a passenger in OD pair $g$ spends on arc $(i,j)$.\\
			$t_{ij}^\text{charge}$ & The charging time that an AEV spends on arc $(i,j)$.\\
			$t_{ij}^\text{drive}$ & The driving time that an AEV spends on arc $(i,j)$.\\
			
			$v$ & Driving speed of AEVs, in kWh/hour.\\
			$\alpha$ & Coefficient for determining number of chargers.\\
			$\eta$ & Charging efficiency of chargers. \\
			$\xi$ & Fuel efficiency of AEVs, in kWh/km.\\
			$\zeta $ & The capital recovery factor that converts the present investment costs into a stream of equal annual payments over the specified lifespan $Y$ at the given discount rate $r$. $\zeta=\frac{r(1+r)^{Y}} {(1+r)_{}^{Y}-1}$.\\
			$\lambda_{g, \tau}^\text{od}$ & Transport demands in OD pair $g$ that depart during hour $\tau$. ${\bm{\lambda}}_{g,\tau}^\text{od}$ is the corresponding vector.\\
			\hline
			\multicolumn{2}{c}{\textbf{Decision variables}}\\
			$x$ & AEV fleet size.\\
			$y_i$ & Number of chargers at node $i$. \\
			$\lambda_{\tau}^\text{drive}$ & Number of AEVs driving on road during hour $\tau$. \\
			$\lambda_{\tau}^\text{park}$ & Total number of AEVs parking during hour $\tau$.\\
			 $\lambda_{i,\tau}^\text{park}$ & Number of AEVs parking at node $i$ during hour $\tau$. $\lambda_{i,\max}^\text{park}$ denotes its upper bound.\\
			$\lambda_{g,\tau}^\text{arr}$ & Number of AEVs in OD pair $g$ that arrive at destination during hour $\tau$.\\
			$\lambda_{g,\tau}^\text{dep}$ & Number of AEVs in OD pair $g$ that depart from origin during hour $\tau$.\\

			$\lambda_{g,ij,\tau}$& Loaded AEVs on arc $(i,j)$ in OD pair $g$ that depart during hour $\tau$. ${\bm{\lambda}}_{g,\tau}$ is the corresponding vector.\\
			$\lambda_{g,ij,\tau}^\text{r}$ &Relocating AEVs on arc $(i,j)$ in OD pair $g$ that depart during hour $\tau$. ${\bm{\lambda}}_{g,\tau}^\text{r}$ is the corresponding vector. \\
			$\lambda_{g, \tau}^\text{odr}$ & Total relocating AEVs in OD pair $g$ that depart during hour $\tau$. ${\bm{\lambda}}_{g,\tau}^\text{odr}$ is the corresponding vector.\\
			
			$\lambda_{g,q,\tau}^\text{path}$ &Loaded AEVs on path $q$ in OD pair $g$ that departs during hour $\tau$.  ${\bm{\lambda}}_{g,\tau}^\text{path}$ is the corresponding vector.\\
			$\lambda_{g,q,\tau}^\text{pathr}$ &Relocating AEVs on path $q$ in OD pair $g$ that departs during hour $\tau$.  ${\bm{\lambda}}_{g,\tau}^\text{pathr}$ is the corresponding vector.\\
			\hline
		\end{tabular}
	\end{small}
\end{table}

\subsection{Objective} 
The objective of the model is to minimize the total costs including the investment for AEV fleet and charging stations, the operation costs for both time, electricity and maintenance.

\subsubsection{Investment Costs for AEVs \& Charging Systems}
The investment costs for purchasing AEVs and installing charging stations can be calculated as follows:
\begin{align}
&c^\text{ev}x + \sum_{i\in \mathcal{I}}  c^\text{sp} y_i,
\end{align}
where, $x$ is the fleet size and $y_i$ is the number of chargers at node $i$. Symbols $c^\text{ev}$ and  $c^\text{sp}$ are the per-unit cost for buying one AEV and installing one charger, in \$, respectively.

\subsubsection{Operation Costs for Time, Electricity and Maintenance}
The operation costs of the AEV system are mainly composed by passenger time, electricity consumption, and maintenance.
The total time that a loaded AEV spends on one arc $(i,j)$ includes two parts: 1) the driving time, $t_{ij}^\text{drive}$, and 2) the charging time at the station, $t_{ij}^\text{charge}$. In terms of costs, only the passenger's time is valuable, while an AEV without a passenger has no time cost. 
We adopt $t_{g,ij}$ to denote the time that a passenger spends on arc $(i,j)$ in OD pair $g$. Then, it is:
\begin{align}
t_{g,ij} &=t_{ij}^\text{drive} +(1-1_{d_g=j})t_{ij}^\text{charge}, \notag\\
&=\frac{l_{ij}}{v} + \left(1-1_{d_g=j}\right)\frac{\xi l_{ij}}{\eta p^\text{sp}}, \quad \forall (i,j)\in \hat{\mathcal{A}}, \forall g\in\mathcal{G},\label{eqn_totalarctime}
\end{align}
where, symbol $v$ represents the AEVs' average driving speed, in km/h; $\xi$ is the fuel efficiency, in kWh/km;\footnote{We assume that the fuel efficiency $\xi$ is homogeneous across the transportation network in this paper. In practice,  different arcs may have different fuel efficiency depending on road conditions. It is trivial to extend our model to consider this fact. Hence, we omitted it for brevity.} $\eta$ is the charging efficiency; $p^\text{sp}$ is the rated charging power, in kW. Hence, $\xi l_{ij}$ is the total electricity that an AEV consumes on arc $(i,j)$. Symbol $1_{d_g=j}$ represents whether arc $(i,j)$'s end node $j$ is also the AEV's destination $d_g$: $1_{d_g=j}$=1, if it is; $1_{d_g=j}$=0, otherwise. Because when a loaded AEV arrives at its destination, it will drop the passenger first before getting charged. Hence, a passenger does not spend any charging time after arriving at its destination (denoted by the last term in (\ref{eqn_totalarctime})).

We assume that the time, electricity and maintenance costs are, respectively, proportional to the time a passenger spent in the vehicle, electricity consumption and vehicle mileage. Hence, the operation costs on arc $(i,j)$ of an AEV driving between OD pair $g$ can be calculated as follows:
\begin{align}
c_{g,ij}^\text{oper}&=c^\text{t} t_{g,ij} + c_{j}^\text{e} \frac{\xi l_{ij}}{\eta} + c^\text{m} l_{ij},\quad\forall (i,j)\in\hat{\mathcal{A}}, \forall g\in\mathcal{G},
\end{align}
where, symbols $c^\text{t}$, $c_j^\text{e}$ and $c^\text{m}$ are the per-unit passenger time cost, in \$/hour, per-unit electricity cost, in \$/kWh, and per-unit maintenance cost, in \$/km, respectively. Note that the electricity costs are location dependent (affected by the locational marginal electricity prices \cite{li2014distribution}) due to unbalanced power supply capabilities. 

\subsubsection{Total Objective}
In summary, the annual expected total investment and operation cost can be formulated as follows:
\begin{align}
&Obj=\min \quad\zeta c^\text{ev}x + \zeta\sum_{i\in \mathcal{I}} c^\text{sp} y_i \notag \\
&+ 365\sum_{g\in \mathcal{G}} \sum_{(i,j)\in\hat{\mathcal{A}}} \sum_{\tau}\left( c^\text{t} t_{g,ij}+ c_{j}^\text{e} \frac{\xi l_{ij}}{\eta} + c^\text{m} l_{ij}\right)\lambda_{g,ij, \tau},\notag \\
&+ 365\sum_{g\in \mathcal{G}} \sum_{(i,j)\in\hat{\mathcal{A}}} \sum_{\tau} \left(c_{j}^\text{e} \frac{\xi l_{ij}}{\eta} + c^\text{m}l_{ij}\right) \lambda_{g,ij, \tau}^\text{r},\label{obj_fleetsize_1}
\end{align}
where, $\zeta $ is the capital recovery factor, which converts the present investment costs into a stream of equal annual payments over the specified lifespan at the given discount rate; $\lambda_{g,ij,\tau}$ and $\lambda_{g,ij,\tau}^\text{r}$ are, respectively, the fraction of loaded and relocating AEVs on arc $(i,j)$ in OD pair $g$ that depart during hour $\tau$. In the objective, the first and second terms are investment costs for the AEV fleet and charging stations, respectively; the third and fourth terms are the operation costs for the loaded and relocating AEVs, respectively. The former has time costs while the latter does not.\footnote{In practice, loaded and relocating AEVs may also have different fuel efficiency. It is trivial to extend our current model to consider this fact.}%(especially for goods transportation when the goods in an loaded AEV can be quite heavy)
\subsection{Constraints}
This sub-section describes the constraints of the planning model. In summary, the dispatch of loaded and relocating AEVs shall be able to satisfy the transport demands and subject to the physical transportation network constraints described in Section \ref{sec_transportation}. To better evaluate the demands for AEVs, we adopt constraints to describe the dynamics of AEV departure, arriving, and parking. In addition, the total fleet size shall be no less than the total number of AEVs driving or parking on the network during any hour, and the installed chargers shall be able to satisfy the peak charging demands at each location.

\subsubsection{AEV Traffic Flow Continuity}
The distribution of the AEV traffic flow should satisfy the transportation network constraints indicated by the expanded network model (\ref{eqn_TS1})--(\ref{eqn_TS2}). In the proposed model, the total AEVs on each path includes the loaded AEVs to satisfy the OD transport demand and the relocating ones. Hence, the expanded network model can be represented as follows:
\begin{align}
&\bm{A}{{\bm{\lambda}}}_{g,\tau}={\bm{\lambda}}_{g,\tau}^\text{od},\qquad\forall g\in \mathcal{G}, \forall \tau, \label{eqn_TS1_2}\\
&\bm{A}{{\bm{\lambda}}}_{g,\tau}^\text{r}={\bm{\lambda}}_{g,\tau}^\text{odr},\qquad\forall g\in \mathcal{G}, \forall \tau, \label{eqn_TS2_2}\\
&{{\bm{\lambda}}}_{g,\tau}\geq 0,~ {{\bm{\lambda}}}_{g,\tau}^\text{r}\geq 0,~ {\bm{\lambda}}_{g,\tau}^\text{odr} \geq 0, \qquad\forall g\in \mathcal{G}, \forall \tau.
\end{align}
where, ${\bm{\lambda}}_{g,\tau}^\text{od}$ and ${\bm{\lambda}}_{g,\tau}^\text{odr}$ are the vectors of transport demands and relocating AEVs in OD pair $g$ departing in hour $\tau$, respectively. Symbol ${\bm{\lambda}}_{g,\tau}$ is the vector of loaded AEV traffic flow on the arcs of the expanded network; its elements are $\lambda_{g,ij,\tau}, \forall (i,j)\in\hat{\mathcal{A}}$; ${\bm{\lambda}}_{g,\tau}^\text{r}$ is the corresponding vector for relocating AEVs whose elements are $\lambda_{g,ij,\tau}^\text{r}, \forall (i,j)\in\hat{\mathcal{A}}$.
Note that the transport demands ${\bm{\lambda}}_{g,\tau}^\text{od}$ are given, while the relocating flow ${\bm{\lambda}}_{g,\tau}^\text{odr}$ are decision variables.

\subsubsection{Relationship between AEVs on Arcs and on Paths}
Assuming that the loaded or relocating AEVs corresponding to demands $(o_g,d_g,\lambda_{g,\tau}^\text{od})$ will only choose the paths in the set $\mathcal{Q}_g$ or $\mathcal{Q}_g^\text{r}$, $\forall g\in\mathcal{G}$, respectively. Then, the relationship between AEVs driving on paths and on arcs are:
\begin{align}
&{{\bm{\lambda}}}_{g,\tau} =\bm{B}_g{\bm{\lambda}}_{g,\tau}^\text{path},\qquad \forall g\in \mathcal{G},\label{eqn_modified1}\\
&{{\bm{\lambda}}}_{g,\tau}^\text{r} =\bm{B}_g^\text{r}{\bm{\lambda}}_{g,\tau}^\text{pathr},\qquad \forall g\in \mathcal{G},\label{eqn_modified2}\\
&{\bm{\lambda}}_{g,\tau}^\text{path} \geq 0, ~ {\bm{\lambda}}_{g,\tau}^\text{pathr}\geq 0,\qquad \forall g\in \mathcal{G},\label{eqn_modified3}
\end{align}
where, ${\bm{\lambda}}_{g,\tau}^\text{path}= \{\lambda_{g,1,\tau}^\text{path}, \lambda_{g,2,\tau}^\text{path}, ..., \lambda_{g,q,\tau}^\text{path},...\}$ is the vector of loaded AEVs in path set $\mathcal{Q}_g$ corresponding to transport demand $(o_g,d_g,\lambda_{g,\tau}^\text{od})$; $\mathbf{B}_g$ is the fundamental path-arc incidence matrix correlated to OD pair $g$ in  $G(\mathcal{I},\hat{\mathcal{A}})$. ${\bm{\lambda}}_{g,\tau}^\text{pathr}$ and $\mathbf{B}_g^\text{r}$ are the corresponding vector and matrix for the relocating AEVs. Based on equations (\ref{eqn_modified1})--(\ref{eqn_modified2}), we can substitute decision variables ${{\bm{\lambda}}}_{g,\tau}$ and ${{\bm{\lambda}}}_{g,\tau}^\text{r}$ by ${\bm{\lambda}}_{g,\tau}^\text{path}$ and ${\bm{\lambda}}_{g,\tau}^\text{pathr}$. If we are able to find small scale path sets, $\mathcal{Q}_g$ or $\mathcal{Q}_g^\text{r}$, $\forall g\in\mathcal{G}$, that cover the paths that AEVs will adopt, the scale of the decision variables can be significantly reduced. Section \ref{sec_algorithm} will discuss how to reduce the complexity of the problem by exploiting the structure of the problem based on these relationships.

\subsubsection{AEV Departure \& Arriving Time} 
As discussed in Section \ref{sec_transportation}, an AEV's total driving and charging time on any path $q$, i.e., $\tau_q$, can be easily calculated \textit{a priori} based on arc length, average driving and charging speed. Hence, there is a unique mapping
 between AEVs departing from the origin, $\lambda_{g,\tau}^\text{dep}$, and AEVs arriving at destination, $ \lambda_{g,\tau}^\text{arr}$, between different hours on each path. This can be described as follows:
\begin{align}
& \lambda_{g,\tau}^\text{dep} = \sum_{q\in\mathcal{Q}_g}\lambda_{g,q,\tau}^\text{path}+ \sum_{q\in\mathcal{Q}_g^\text{r}}\lambda_{g,q,\tau}^\text{pathr}, \qquad \forall g\in\mathcal{G}, \forall \tau,\\
& \lambda_{g,\tau}^\text{arr} = \sum_{q\in\mathcal{Q}_g}\lambda_{g,q,\tau-\tau_q}^\text{path}+\sum_{q\in\mathcal{Q}_g^\text{r}}\lambda_{g,q,\tau-\tau_q}^\text{pathr}, \quad \forall g\in\mathcal{G}, \forall \tau.
\end{align}
\subsubsection{AEV Driving on Road}
Based on the above analysis on AEV departure and arriving time on paths, it is straightforward to calculate the number of AEVs on road, as follows:
\begin{align}
\lambda_{\tau}^\text{drive} = &\sum_{\tau_0}\sum_{g\in \mathcal{G}} \bigg(\sum_{q\in\mathcal{Q}_g}1_{\tau\in[\tau_0,\tau_0+\tau_q]} \lambda_{g,q,\tau_0}^\text{path} +\notag\\
& \sum_{q\in\mathcal{Q}_g^\text{r}}1_{\tau\in[\tau_0,\tau_0+\tau_q]} \lambda_{g,q,\tau_0}^\text{pathr}\bigg), \qquad \forall \tau,
\end{align}
where, $1_{\tau\in[\tau_0,\tau_0+\tau_q]} =1$ if the AEVs on path $q$ are driving between time $[\tau_0,\tau_0+\tau_q]$; $1_{\tau\in[\tau_0,\tau_0+\tau_q]} =0$, otherwise. 

\subsubsection{AEV Parking Dynamics}
In practice, the AEV fleet operator may initially locate a number of AEVs at each location for future demands. During the operation, the number of vehicles parking at each node will change. The corresponding dynamics can be represented as follows:
\begin{align}
& \lambda_{i,\tau}^\text{park} = \lambda_{i,\tau-1}^\text{park} + \sum_{g\in \mathcal{G}} \left({1_{d_g=i}\lambda_{g,\tau}^\text{arr}} - { 1_{o_g=i}\lambda_{g,\tau}^\text{dep}}\right) ,  \forall i\in\mathcal{I},\forall \tau, \label{eqn_park1}\\
& \lambda_{i,T}^\text{park} \geq \lambda_{i,0}^\text{park},\qquad \forall i\in\mathcal{I},\label{eqn_park2}\\
& \lambda_{i,\tau}^\text{park} \leq \lambda_{i,\max}^\text{park},\qquad \forall \tau, \forall i\in\mathcal{I},\label{eqn_park3}\\
&\lambda_{\tau}^\text{park} = \sum_{i\in\mathcal{I}}\lambda_{i,\tau}^\text{park}, \qquad\forall\tau,\label{eqn_park4}
\end{align}
where, symbol $\lambda_{i,\tau}^\text{park}$ denotes the number of AEVs parking at node $i$ during hour $\tau$; $\lambda_{i,0}^\text{park}$ represents the initial number of AEVs that are located at node $i$, which is also a decision variable for the system operator. Equations $\sum_{g\in \mathcal{G}} {1_{d_g=i}\lambda_{g,\tau}^\text{arr}}$ and $\sum_{g\in \mathcal{G}} {1_{d_g=i}\lambda_{g,\tau}^\text{dep}}$ calculate the AEVs that arrive at and depart from node $i$ during hour $\tau$, respectively. Equation (\ref{eqn_park1}) determines the relationship between numbers of AEV parking, arrival and departure. To avoid myopic fleet relocating, we use equation (\ref{eqn_park2}) to require that the number of parking AEVs at each node will return to the corresponding initial value after the operation period, $T$, e.g., one day. Equation (\ref{eqn_park3}) constrains the maximum number of AEVs parking at each location due to parking space limitation, in which symbol $\lambda_{i,\max}^\text{park}$ denotes the upper bound. Symbol $\lambda_{\tau}^\text{park}$ represents the total number of parking AEVs during hour $\tau$, which is calculated by (\ref{eqn_park4}).

\subsubsection{AEV Fleet Size}
The total number of AEVs should be higher than the summation of those parking and driving:
\begin{align}
&x \geq \lambda_{\tau}^\text{drive} + \lambda_{\tau}^\text{park}, \qquad \forall \tau.
\end{align}
%\footnote{In this constraint, the number of AEV fleet is only designed for the average system demands or assuming the demands are not stochastic. In practice, we may need to invest in more AEVs to guarantee the system's quality of service. In that scenario, the service level model can be adopted.}

\subsubsection{AEV Charging Station}
At each node $j$, there should be enough chargers to guarantee adequate quality of service:
\begin{align}
&y_j \geq \alpha\sum_{g\in \mathcal{G}} \sum_{\{i|(i,j)\in \hat{\mathcal{A}}\}}\frac{\xi l_{ij}}{\eta p^\text{sp}}\left(\lambda_{g,ij}+ \lambda_{g,ij}^\text{r}\right), ~ \forall j\in {\mathcal{I}},\label{eqn_fleetsize_QoS}
\end{align}
where, $\frac{\xi l_{ij}}{\eta p^\text{sp}}$ is the required charging time of an AEV on arc $(i,j)$. Symbols $\lambda_{g,ij}$ and $\lambda_{g,ij}^\text{r}$ represent the loaded and relocating AEV flow on arc $(i,j)$ during the peak hour, respectively. Symbol $\alpha$ is a coefficient that is higher than 100\% to make the planning slightly conservative.\footnote{When $\alpha$=100\%, $y_i$ is the number of chargers to satisfy the mean demands.} 

\subsection{Complete Fleet Sizing \& Charging System Planning Model}

The above formulations form the joint AEV fleet sizing and charging station planning model summarized as follows:
\begin{align}
\textbf{P1:}~~&\min \quad \text{(\ref{obj_fleetsize_1})} \qquad \text{s.t.:\quad (\ref{eqn_TS1_2})--(\ref{eqn_fleetsize_QoS})}.\notag
\end{align}
It is a mixed-integer linear program. 

For goods transportation, neither loaded nor relocating AEVs have passengers or drivers, hence, there will be no time cost for them. The objective in (\ref{obj_fleetsize_1}) is modified as:
\begin{align}
&Obj=\min \quad\zeta c^\text{ev}x + \zeta\sum_{i\in \mathcal{I}} c^\text{sp} y_i \label{obj_fleetsize_2} \\
&+ 365\sum_{g\in \mathcal{G}} \sum_{(i,j)\in\mathcal{A}} \sum_{\tau}\left(c_{j}^\text{e} \frac{\xi }{\eta} + c^\text{m}\right)l_{ij}\left(\lambda_{g,ij, \tau}+\lambda_{g,ij, \tau}^\text{r}\right).\notag
\end{align}

Note that without time cost in goods transportation, an AEV may detour more in order to get charged at lower electricity prices or allow the system to install less chargers. However, long driving and charging time will still affect an AEV's utilization so that it will not detour unrestrained.

\section{Solution Approach}\label{sec_algorithm}
Because an AEV may theoretically choose any path available to travel between any OD pair in the network, the proposed model \textbf{P1} can be computationally expensive when the network is large. 
However, in practice, the choices of an AEV between one OD pair can be quite limited because detouring from the shortest path to other ones could be expensive due to high time, electricity and maintenance costs. Therefore, if we are able to find the smallest path sets $\mathcal{Q}_g$ and $\mathcal{Q}_g^\text{r}$, $\forall g\in \mathcal{G}$, in equations (\ref{eqn_modified1})--(\ref{eqn_modified2}) that contain all the paths that will be used by AEVs, we can reduce its complexity significantly without sacrificing any optimality.

Finding these path sets is challenging. Nevertheless, in this section, we propose some supersets of $\mathcal{Q}_g$ and $\mathcal{Q}_g^\text{r}$, $\forall g\in \mathcal{G}$, by exploiting the structure of this problem and deriving valid \textit{cuts}, i.e., inequality constraints, to eliminate useless paths between each OD pair. Utilizing the proposed supersets can still significantly reduce the computational efficiency of the problem \textbf{P1} without sacrificing optimality. Furthermore, we also propose subsets of the aforementioned supersets to approximately solve the problem by tuning a confidence factor of the approximation to balance computational efficiency and optimality. As a result, the problem can be still tractable in large-scale networks with guaranteed solution quality.

\subsection{Cuts for Loaded AEVs}\label{cuts_1}
\subsubsection{Warm-up Solution}
First, we ignore the investment costs for AEV fleet and charging systems or assume that they are negligible compared to operation costs. Then, an AEV will only choose the path that minimizes its operation costs between an OD pair, which is referred to as the \textit{min-operation path} hereinafter. As a result, each OD pair only corresponds to one path so that the planning problem will reduce to a small scale problem that can be efficiently solved by an off-the-shelf solver. However, the result may be sub-optimal when operation costs are not dominating the system or the company has limited budget and hopes to save its investment costs. 

We will adopt the above solution with min-operation paths as a warm-up based on which we will exploit the structure of the problem and propose valid cuts for it. 

\subsubsection{Cuts Relevant to Fleet Size}
For the fleet sizing problem, when an AEV is driving or charging, it is not available for relocating or servicing other demands. Hence, any time delay caused by detour may lead to opportunity cost for investing in a larger fleet. Therefore, it is possible to reduce the system's total costs by routing some AEVs to drive on paths with shorter time rather than min-operation paths. 

We use $q$ to denote the index of paths between one OD pair $g$, $q\in\mathcal{Q}_g$; $0$ to denote the index of the min-operation path. Let $\hat{\mathcal{A}}_q$ denote the set of the arcs on path $q$. Then, if a fleet of AEVs, $\lambda$, drive on path $q$, the incurred operation costs are:
\begin{align}
&C_{q} = 365\sum_{(i,j)\in\hat{\mathcal{A}}_q} \left( c^\text{t} t_{g,ij}+ c_{j}^\text{e} \frac{\xi l_{ij}}{\eta}  + c^\text{m}\right)\lambda.
\end{align}
Note that for goods transportation, $c^\text{t}=0$. The total driving and charging time on path $q$ is $t_q=\sum_{(i,j)\in\hat{\mathcal{A}}_q} t_{g,ij}$. 

Based on the above analysis, when we ignore the investment costs in charging systems, we have the following proposition:

\textit{Proposition 1: } The maximum total system cost reduction for \textbf{P1} when $c^\text{sp}=0$ by detouring $\lambda$ AEVs from min-operation path $0$ to any other path $\forall q\in {\mathcal{Q}\setminus 0}$ is upper bounded by:
\begin{align}
&\Delta C_{q,1}=
\left\{
\begin{array}{l}
C_{0} - C_{q},\qquad\qquad ~\text{if}~t_{q} \geq  t_{0}\\
C_{0} - C_{q} + \zeta c^\text{ev}\lambda,~~~\text{if}~t_{q} <  t_{0}\end{array}\right..
\end{align}

\textit{Proof:} The proof for this proposition is intuitive: 1) When $t_{q} \geq  t_{0}$, detouring one AEV from path $0$ to path $q$ will not reduce the fleet size because the AEV will arrive at its destination later than before. Hence, the total system cost will at least increase by $C_q-C_0\geq 0$. 2) When $t_{q} <  t_{0}$, $\lambda$ AEVs will arrive at their destination and be able to service other demands earlier. The most optimistic scenario is that the earlier arrived AEVs happen to satisfy the gap between demands and supply at that location during that specific time interval so that $\lambda$ AEVs can be saved. Hence, the maximum possible cost reduction is the difference between the saved fleet investment $\zeta c^\text{ev}\lambda$ and the operation cost increment due to detouring $C_q-C_0$, i.e., $C_{0} - C_{q} + \zeta c^\text{ev}\lambda$.

\begin{remark}
	Based on Proposition 1, an AEV may only choose a path $q$ when $\Delta C_{q,1}>0$.  Hence, all other paths can be eliminated, which constructs a set of valid cuts for \textbf{P1}. 
\end{remark}

\subsubsection{Cuts Relevant to Charging Systems}
Similarly, we can also derive cuts considering charging system investments. 
If we ignore the fleet investment costs, we have:

\textit{Proposition 2: } The maximum total system cost reduction for \textbf{P1} when $c^\text{ev}=0$  by detouring $\lambda$ AEVs from min-operation path $0$ to any other path $\forall q\in {\mathcal{Q}\setminus 0}$ is upper bounded by:
\begin{align}
&\Delta C_{q,2}=C_{0} - C_{q} + \zeta c^\text{sp}\alpha \sum_{(i,j)\in\hat{\mathcal{A}}_0^{*}} \frac{\xi l_{ij}}{\eta p^\text{sp}}\lambda,
\end{align}
where, set $\hat{\mathcal{A}}_0^{*} = \hat{\mathcal{A}}_0\setminus (\hat{\mathcal{A}}_0 \cap \hat{\mathcal{A}}_q)$ represents the arc set on path $0$ that the AEVs will no longer visit after detouring so that the corresponding charging demands will be satisfied elsewhere.

\textit{Proof:} The proof for Proposition 2 is also intuitive: after detouring, only those chargers that are installed on arcs in $\hat{\mathcal{A}}_0^{*}$ that are no longer visited by fleet, $\lambda$, may be saved. 

\begin{remark}
	Based on Proposition 2, an AEV may only choose path  $q$ when $\Delta C_{q,2}>0$. Thus, all other paths can be eliminated. This also constructs valid cuts for model \textbf{P1}. 
\end{remark}

\subsubsection{Cuts Relevant to Both Fleet Size \& Charging Systems}
Based on the above analysis, it is straightforward that:

\textit{Proposition 3:} The maximum total system cost reduction for \textbf{P1} by detouring $\lambda$ AEVs from min-operation path $0$ to any other path $\forall q\in {\mathcal{Q}\setminus 0}$ is upper bounded by:
\begin{align}
&\Delta C_q =
\left\{
\begin{array}{l}
C_{0} - C_{q} + d,\qquad\qquad ~\text{if}~t_{q} \geq  t_{0}\\
C_{0} - C_{q} + \zeta c^\text{ev}\lambda + d,~~~\text{if}~t_{q} <  t_{0}\end{array}\right.,
\end{align}
where, $d = \zeta c^\text{sp}\alpha \sum_{(i,j)\in\hat{\mathcal{A}}_0^{*}} \frac{\xi l_{ij}}{\eta p^\text{sp}}\lambda$.

In practice, we can first adopt a \textit{$k$ shortest path routing algorithm} \cite{yen1971finding} to identify $k$ paths with the least operation costs for each OD pair; then, remove those paths with $\Delta C_q \leq 0$ based on Proposition 3 to construct a superset of $\mathcal{Q}_g$, $\forall g\in \mathcal{G}$. As a result, we need only consider these super path sets in the optimization model so that the problem scale can be significantly reduced and optimality is still retained.

If the supersets are still too huge due to the network's large scale, we can remove all the paths that have $\Delta C_q \leq gap\times obj_0$, where $obj_0$ is the total objective value of the warm-up solution, and $gap$ is a confidence factor. Apparently, by letting $gap>0$, we are able to reduce the number of paths remaining in the path sets, which may result in a sub-optimal solution. By adjusting $gap$'s value, we can balance computational efficiency with optimality of the problem effectively.

\subsection{Cuts for Relocating AEVs}\label{cuts_2}
For relocating, an AEV parked at one node may be relocated to any other node in the transportation network, and one AEV relocating from one origin to a destination may also choose different paths. As a result, the routing of relocating AEVs can be as complex as the loaded ones. However, based on the assumption [A\ref{as_rebalance}], we have the following proposition:

\textit{Proposition 4: } Relocating an AEV on path $q$ with arc set $\hat{\mathcal{A}}_q$ is equivalent to relocating one AEV on each arc (from the start node to the end node of the arc) in $\hat{\mathcal{A}}_q$ successively. 

\textit{Proof:} The proof for Proposition 4 is also intuitive. We take a three-arcs path $(ijk)$ for example, assume that we need to allocate one AEV from $i$ to $k$ on path $(ijk)$. According to assumption [A\ref{as_rebalance}], the AEV will depart with a full battery and get fully charged after arriving at node $j$ and node $k$, respectively. It is obvious that, if we first relocate the AEV from $i$ to $j$ and then relocate it again from $j$ to $k$, the AEV will spend the same length of time (for driving and charging) with the case when we relocate it from $i$ to $k$ on path $(ijk)$.

\begin{remark}
	Based on Proposition 4, we can require that the operator only relocate AEVs  between adjacent nodes to significantly reduce computational complexity without losing any optimality, which constructs the cuts for relocating AEVs. As a result, the usable path set $\mathcal{Q}_g^\text{r} = \{(o_g,d_g)\}$ if $(o_g,d_g)\in \hat{\mathcal{A}}$; $\mathcal{Q}_g^\text{r} = \emptyset$, otherwise, $\forall g\in \mathcal{G}$. 
\end{remark}

Because the scale of the AEV routing and relocating decision variables is proportional to number of usable paths of OD pairs, by adding the above cuts to eliminate useless cuts will help significantly reduce the complexity of the proposed model. As a result, the model can be efficiently solved by the branch-and-bound algorithm.

\section{Case Studies}\label{sec_case}
\begin{figure}
	\centering
%	\vspace{-4mm}
	\includegraphics[width=0.5\columnwidth]{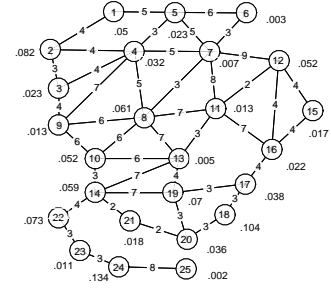}
	\vspace{-3mm}
	\caption{A 25-node transportation network used for the case study. The number in each circle is the node index. The number on each arc represents the distance between the corresponding two nodes and the per-unit distance is 10 km. The decimal next to each node is its traffic flow gravitation\cite{Transportation25_Simchi1988}.}
	\label{fig_trans25}
	\vspace{-5mm}
\end{figure}

This section considers a 25-node  transportation network (see Fig.~\ref{fig_trans25}) to illustrate the proposed planning method. The gravity spatial interaction model was used to generate a time varying OD demands based on node weights and arc distances\cite{Transportation25_Simchi1988}. The details of the trip demand information can be found in our previous publication \cite{Hongcai_SOCP_Zhang2017}. We assume the transport demands on this network are about 15,000 trips/day, with 2,250 trips/hour during the peak hour. 

\subsection{Parameter Settings}
We utilize the AEV and charger cost parameters published by Bauer et. al.\cite{Bauer2018}. The per unit cost to purchase an AEV is $c^\text{ev}$=30,000+200$B$, in \$, where $B$ is the AEV's battery capacity, in kWh; The per unit cost to install a charger is $c^\text{sp}$=(700 + 15Y)$p^\text{sp}$, where $Y$ is the life time of a charger, in year, $p^\text{sp}$ is the rated charging power, in kW. The vehicle efficiency $\xi$=0.155+0.00037$B$, in kWh/km. The maintenance cost $c^\text{m}$=0.025 kWh/km. We also assume the per-unit time cost for a passenger is  $c^\text{t}$=22.62 \$/hour, which is the average hourly earnings of private-sector production and nonsupervisory employees in the US\cite{USsalary}. The charging efficiency $\eta$=0.92, the capital recovery factor $\zeta=\frac{r(1+r)^{Y}} {(1+r)_{}^{Y}-1}$, where $Y$=15 year, discount rate $r$=0.08\cite{Hongcai_SOCP_Zhang2017}.
We assume the rated charging power of a charger is $p^\text{sp}$=100 kW and the average driving speed is 100 km/h. 
Because of congestion in power networks, the electricity supply costs in different areas may also be different\cite{li2014distribution}. Hence, we assume heterogeneous electricity prices across the transportation network, i.e., 0.12 \$/kWh at nodes 6-8, 11-13, 16; 0.20 \$/kWh at nodes 10, 14, 21-25, and 0.30 \$/kWh at other nodes.

For the base case, we assume the AEVs' battery capacities are all 75 kWh, which is equal to that of a TESLA Model 3. However, we will not model the driving range constraints of the AEVs described in Section \ref{sec_transportation} based on this value because of four major reasons: 1) First, an AEV should conserve sufficient residual battery electricity to make sure that it can reach to the next charging station safely; 2) Second, a battery's charging speed will dramatically slow down when its state-of-charge gets high so that an AEV will not tend to get fully charged to save time; 3) Thirdly, the fuel efficiency of an AEV may change with traffic or environmental conditions; and 4) Lastly, the battery capacity of an AEV may degrade during the lifetime.
Therefore, we only use 60 kWh battery capacity to calculate the driving range of AEVs in the transportation network. For other scenarios with different battery capacities, we also deduct 15 kWh from the nameplate capacity to calculate driving range because of the aforementioned reasons.

The expanded transportation network has 25 nodes with 590 arcs and 25$\times$24=600 OD pairs. The scales of major decision variables $\lambda_{g,ij,\tau}$ and $\lambda_{g,ij,\tau}^\text{r}$ are both in the magnitude of 8 million. We first substitute them by $\lambda_{g,q,\tau}^\text{path}$ and $\lambda_{g,q,\tau}^\text{pathr}$ using equations (\ref{eqn_modified1})--(\ref{eqn_modified2}). Then, we set the maximum path number of one OD pair $k$=150 and $gap$=$10^{-4}$ and adopt the cuts in Section \ref{sec_algorithm} to eliminate the redundant paths. As a result, by introducing the cuts in Section \ref{cuts_1}, the scale of decision variables $\lambda_{g,q,\tau}^\text{path}$ is only in the magnitude of 0.1 million. By introducing the cuts in Section \ref{cuts_2}, the number of  $\lambda_{g,q,\tau}^\text{pathr}$ is only 590$\times$24=14160. In summary, by utilizing the proposed cuts, the scale of decision variables is reduced by about 98\%. The original problem is intractable, but the new problem with the proposed cuts can be solved by Gurobi on a desktop with a 36-core Intel Xeon Gold 6140 CPU and 64 GB memory in less than 30 minutes.

\begin{table*}
	\small
	\centering
	\vspace{-4mm}
	\begin{small}
		\caption{Optimization results in different cases}
		\vspace{-2mm}
		\begin{tabular}{cccc|cccccc}
			\hline
			\multirow{2}*{Case}&\multirow{2}*{Strategy}&\multicolumn{2}{c|}{Planning result}&\multicolumn{6}{c}{Equivalent annual costs (M\$)}\\
			\cline{3-10}
			&&Fleet size& No. of chargers & Investment & Driving time& Charging time & Electricity & Maintenance &Total\\
			\hline
			\parbox[t]{2mm}{\multirow{4}{*}{\rotatebox[origin=c]{90}{Passenger}}} 
			&MinTime& 2512 & 390 & 2.04 & 102.0 & 0.00 & 20.83 & 12.74 & 137.58 \\ 
			﻿&MinOperation& 2512 & 390 & 2.04 & 102.0 & 0.00 & 20.83 & 12.74 & 137.58 \\  
			﻿&NoRelocation& 3162 & 628 & 2.74 & 102.0 & 0.02 & 21.02 & 12.72 & 138.47 \\  
			﻿&Optimal& 2497 & 330 & 1.95 & 102.0 & 0.02 & 20.82 & 12.74 & 137.51 \\ 
			\hline
			\parbox[t]{2mm}{\multirow{4}{*}{\rotatebox[origin=c]{90}{Goods}}} 
			&MinTime& 2512 & 390 & 2.04 &  0.0 & 0.00 & 20.83 & 12.74 & 35.60 \\  
			﻿&MinOperation& 2540 & 397 & 2.06 &  0.0 & 0.00 & 18.15 & 12.91 & 33.12 \\  
			﻿&NoRelocation& 3167 & 551 & 2.64 &  0.0 & 0.00 & 18.35 & 12.88 & 33.87 \\ 
			﻿&Optimal& 2505 & 252 & 1.86 &  0.0 & 0.00 & 18.18 & 12.89 & 32.93 \\  
			\hline
		\end{tabular}\label{tab_result}
	\end{small}
	\vspace{-4mm}
\end{table*}

\begin{figure*}
	\centering
%	\vspace{-4mm}
	\includegraphics[width=2\columnwidth]{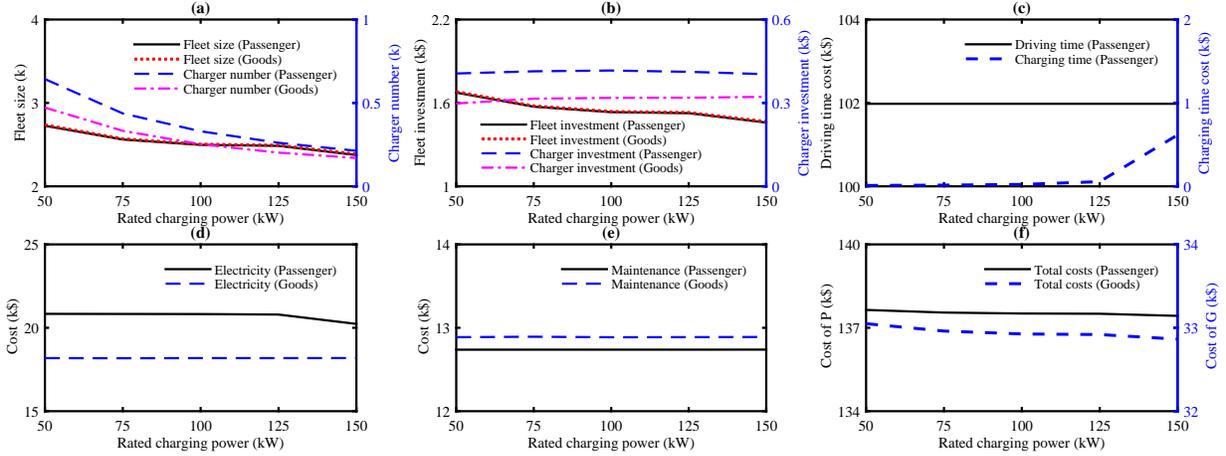}
	\vspace{-4mm}
	\caption{Planning results with different rated charging power. (a), sizes of the EV fleet and charging systems. (b), investment costs. (c), time costs (Note that there is no time costs for goods transportation), (d) electricity costs, (e) maintenance costs, (f) total costs.}
	\label{fig_power}
	\vspace{-5mm}
\end{figure*}
\subsection{Results and Analysis}
We propose three benchmarks to validate the efficacy of the proposed strategy, as follows: 1) MinTime: the AEVs simply choose their shortest paths (with minimum driving time) to drive; 2) MinOperation: the AEVs simply choose their min-operation paths (with minimum operation costs) to drive; and 3) NoRelocation: the AEV fleet operator does not actively relocate the AEVs during the day; instead, it only relocates the AEVs to their initial states after one day of operation.

The planning results for both passenger and goods transportation in the base case are summarized in Table~\ref{tab_result}. Generally, for the system with 15,000 trips/day, a total of 2,500+ AEVs and 300+ chargers are required. That means one AEV can satisfy 6 demands and one charger can fuel about 8 AEVs on average each day. 

As expected, the proposed strategy is the most beneficial one for both passenger and goods transportation. Compared with the MinTime, MinOperation and NoRelocation Strategies, the proposed strategy can significantly reduce the required investments, especially for chargers, by optimally routing AEVs and scheduling their charging locations. Though the total costs reduction may be marginal (because of high operation costs share), the proposed strategy can reduce investment costs by at least 4.4\% and 8.8\% compared with benchmarks for passenger and goods transportation, respectively. This can help the company to reduce high capital investment that can be vitally important when its budget is limited. 

When adopting the MinTime strategy, AEVs tend to choose the paths that will minimize the total passenger time on the road. However, because it overlooks other operation costs, its electricity costs will be higher. This will deteriorate the optimization results for goods transportation when electricity is the dominant cost factor leading to a 8.1\% increase of the total costs compared with the proposed strategy. However, its planning results are pretty close to the MinOperation strategy. 

The performance of the MinOperation strategy is quite close to the proposed one. The operation costs for both of the two strategies are approximately equal. However, the latter does remarkably advance the former in terms of investment costs. This indicates that the proposed strategy can effectively coordinate the route choices of AEVs in different OD pairs (with heterogeneous driving behaviors) to reduce charging demand peaks without sacrificing operation costs too much. 

Relocating AEVs during the operation is quite important for promoting AEV utilization and reduce system costs. Without actively relocating AEVs, the number of AEV fleet could be 26.6\% and 26.4\% higher for passenger and goods transportation, respectively. At the same time, the number of chargers could be doubled for both cases so that the investment costs can be 40.5\% and 41.9\% higher, respectively. 

Generally, the proposed strategy is more effective for goods transportation with less operation costs share. With the same quantity of mobility demands, the investment and electricity costs for goods transportation are 5.1\% and 12.7\% lower than those of passenger transportation, respectively. That is because AEVs are more willing to detour when carrying goods than passengers. For the latter, the passenger time cost is the dominant cost factor so that detouring, which increase driving and charging time, will be expensive. However, for the former, there is no passenger time cost and electricity cost is the dominant factor that is comparatively cheaper. Hence, AEVs are more willing to detour to charge cheaper electricity and this also reduces investments for both fleet with chargers. 

\subsection{Sensitivity Analysis}

\subsubsection{Rated Charging Power}
As mentioned earlier, because charging is much slower than gas filling, it can lead to  considerable downtime (up to 45 minutes or more). With the fast development of battery storage and power electronic technologies, the rated charging power of chargers will also increase. The planning results of the fleet size and charging system as well as the corresponding costs structure with different rated charging power are illustrated in Fig.~\ref{fig_power}.

When the power increases, a same charger can satisfy more demands and the average charging time of AEVs will be lowered. Thus, the demands for AEVs and chargers are both reduced remarkably. Though installing high-power chargers is quite expensive, the experiment results show that investing in higher power chargers is more economical.\footnote{Note that, in these experiments, we did not count the possible demands for power system upgrades to support the chargers.}

\begin{figure*}
	\centering
	\vspace{-4mm}
	\includegraphics[width=2\columnwidth]{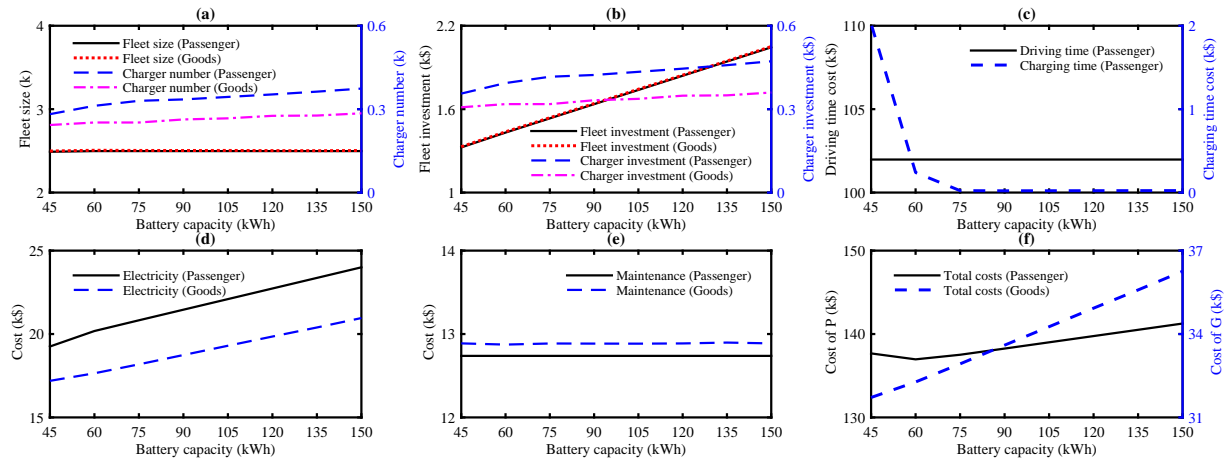}
	\vspace{-5mm}
	\caption{Planning results with different battery capacities. (a), sizes of the EV fleet and charging systems. (b), investment costs. (c), time costs (Note that there is no time costs for goods transportation), (d) electricity costs, (e) maintenance costs, (f) total costs.}
	\label{fig_capacity}
	\vspace{-5mm}
\end{figure*}
\begin{figure*}
	\centering
	%	\vspace{-10mm}
	\includegraphics[width=2\columnwidth]{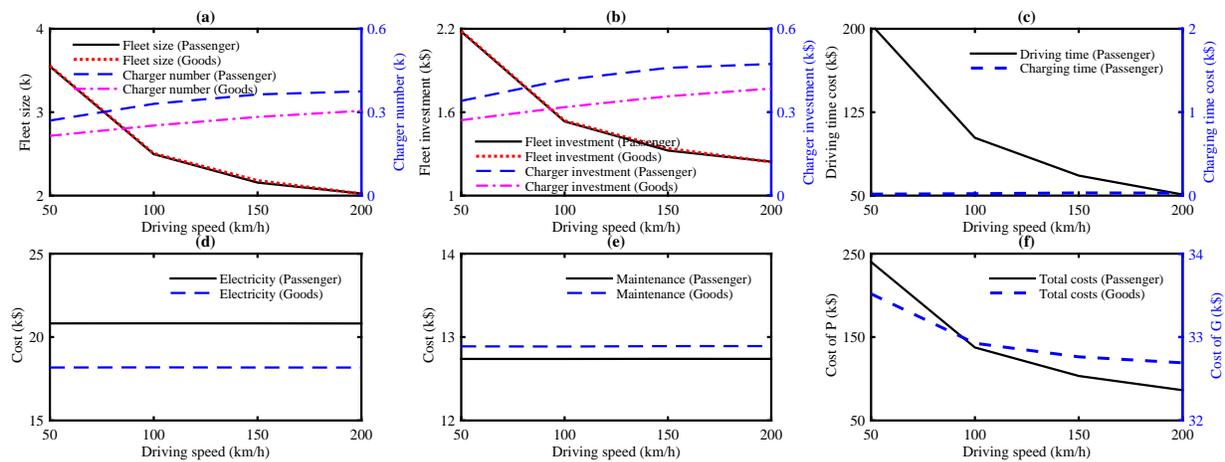}
	\vspace{-5mm}
	\caption{Planning results with different average driving speed. (a), sizes of the EV fleet and charging systems. (b), investment costs. (c), time costs (Note that there is no time costs for goods transportation), (d) electricity costs, (e) maintenance costs, (f) total costs.}
	\label{fig_speed}
	\vspace{-4mm}
\end{figure*}
\subsubsection{Battery Capacity}
Battery capacity is another factor that may significantly affect the system's performance. The summary of the planning results with different battery capacities are illustrated in Fig.~\ref{fig_capacity}. 
When AEVs' battery capacity rises up, their fuel efficiency will reduce due to added mass. As a result, AEVs' electricity consumption will increase so that the demands for chargers and electricity costs will rise up. 
For goods transportation, this will increase the system's total costs and not be economical. 

However, for passenger transportation, AEVs with higher battery capacity will charge less often on road so that they can save passengers' time costs. This cost reduction may counteract the cost increment due to higher electricity consumption. As is shown in Fig.~\ref{fig_capacity}(f), the total costs of the system for passenger transportation first decrease and then increase with the increase of battery capacity. The 60 kWh battery capacity leads to the lowest cost.

\subsubsection{Average Driving Speed}
For AEVs, it is possible that the mobility technology will relax the current transportation regulations, e.g., speed limit. As a result, AEVs may be allowed to drive with higher speed. We validate how this will impact the planning results, which are summarized in Fig.~\ref{fig_speed}. As expected, with higher average driving speed, the demands for fleet size and operation time costs will both be significantly reduced. The total system costs for both passenger and goods transportation can be dramatically reduced.

\section{Conclusion}\label{sec_conclusion}

We propose a planning model to right-size an AEV fleet and charging stations at least cost while meeting transport demands for both passenger and goods transportation. Based on the proposed model we studied various parameters' impact on the system performance. 
Numerical experiments prove that the proposed strategy can effectively balance the investment costs and operation costs of the AEV system. By considering future routing and relocating operations, we can remarkably reduce the investment costs at the planning stage. The planning results also show that AEVs tend to detour more for goods transportation compared with passenger transportation because detouring with passengers may lead to significant time costs. 

Adopting higher power chargers will reduce the downtime of AEVs due to charging, so that the turn-over rate of chargers and the utilization of AEVs are enhanced. Thus, higher charging power results in a lower required number of chargers and AEV fleet size. However, in practice, the company may need to upgrade expensive power supply infrastructure in order to support high power chargers and AEV batteries may age faster using higher power chargers. This may cancel out the mentioned benefits, which needs further study.

Larger battery capacity leads to lower fuel efficiency (which means more electricity consumption) and less charging time with passengers on road. As a result, it is not reasonable to invest in larger AEV batteries for goods transportation because there is not passenger time cost. However, the battery capacity should be carefully designed for passenger transportation balancing fuel consumption increase and charging time reduction.

In this paper, we assume whenever an AEV arrives at its destination, it will get fully charged before driving to another location, which is a mild one for long-distance transportation. However, if the mobility demands are short-distance ones, i.e., taxi-hailing services in urban areas, this assumption will make our planing results conservative. Relaxing this constraint will be our future focus. This paper also assumes that the routing of AEVs will not affect the traffic conditions, this may not be true in practice. In our future works, we will also further study AEVs' potential impact on traffic congestion and the corresponding influence on planning.
In addition, the proposed model is a cost minimization problem and does not account for profit. A net-profit maximizing problem would also need to account for passenger or goods supply and demand economics. This is also a potential future research topic.

\bibliographystyle{ieeetr}
\bibliography{bibliography}

\begin{IEEEbiography}[{\includegraphics[width=1in,height=1.25in,clip,keepaspectratio]{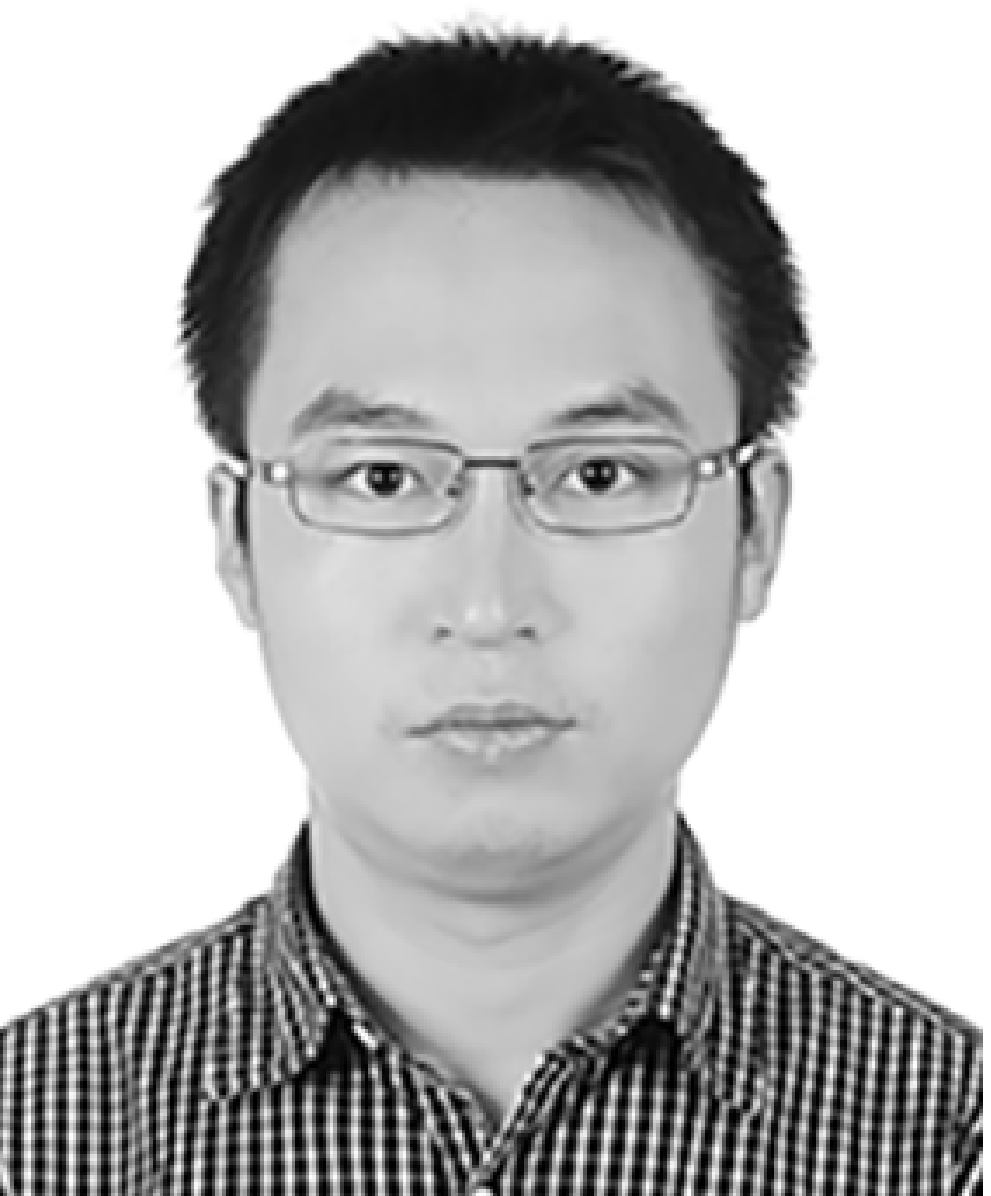}}]{Hongcai Zhang}(S'14--M'18)
	received the B.S.  and Ph.D. degree in electrical engineering from Tsinghua University, Beijing, China, in 2013 and 2018, respectively. In 2016--2017, he worked as a visiting student researcher in the Energy, Controls, and Applications Lab at University of California, Berkeley, where he is currently working as a postdoctoral scholar.
	His current research interests include optimal  operation and optimization of power and transportation systems, grid integration of distributed energy resources. 
\end{IEEEbiography}

\begin{IEEEbiography}[{\includegraphics[width=1in,height=1.25in,clip,keepaspectratio]{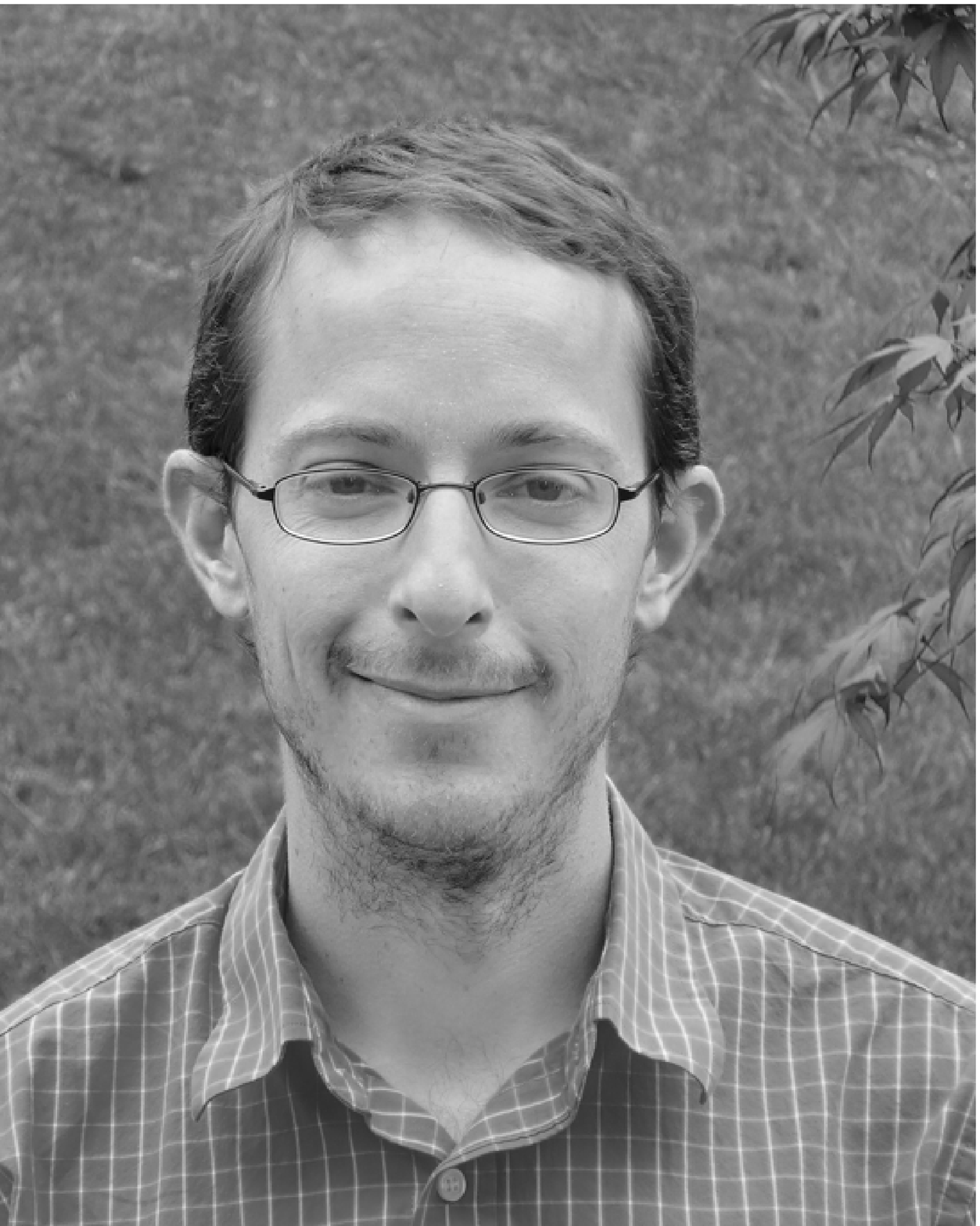}}]{Colin J. R. Sheppard}
	received the B.S. degree in symbolic systems from the Stanford University, Stanford, CA, USA, in 2001, the M.S. degree in environmental resources engineering from the Humboldt State University, Arcata, CA, USA, in 2009, and is currently pursuing the Ph.D. degree in transportation engineering at the University of California, Berkeley, CA, USA. He is a Senior Scientific Engineering Associate with the Sustainable Energy Systems Group at Lawrence Berkeley National Laboratory, Berkeley, CA, USA. He is developing analytical tools to understand the energy and mobility impact of vehicle automation, electrification, and sharing, in addition to the potential of plug-in electric vehicles to provide grid services through an integrated systems approach. 
\end{IEEEbiography}

\begin{IEEEbiography}[{\includegraphics[width=1in,height=1.25in,clip,keepaspectratio]{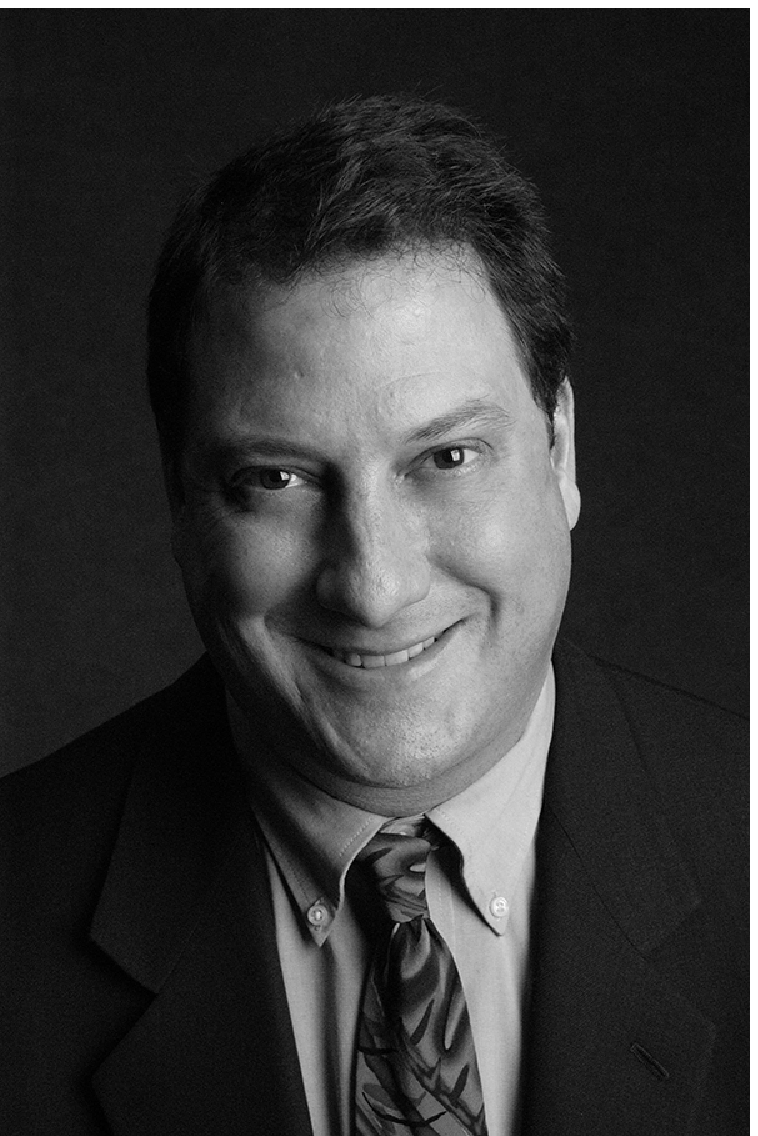}}]{Timothy E. Lipman}
	 is an energy and environmental technology, economics, and policy researcher with the University of California - Berkeley. He is currently serving as Co-Director of the Transportation Sustainability Research Center (TSRC), based at the Institute of Transportation Studies. He is also a Research Affiliate with the Lawrence Berkeley National Laboratory (LBNL). Dr. Lipman’s research focuses on electric vehicles, fuel cell technology, combined heat and power systems, renewable energy, and electricity and hydrogen production and distribution infrastructure. He completed a Ph.D. degree in Environmental Policy Analysis with the Graduate Group in Ecology at UC Davis (1999) and also holds an MS degree from UC Davis in Transportation Technology and Policy (1998) and a BA degree from Stanford University (1990). He is Chair of the Alternative Transportation Fuels and Technologies Committee of the Transportation Research Board of the National Academies of Science and Engineering, a member of the Advisory Committee for the Bay Area Air Quality Management District, and on the editorial boards of Transportation Research-D and the International Journal of Sustainable Engineering.
\end{IEEEbiography}

\begin{IEEEbiography}[{\includegraphics[width=1in,height=1.25in,clip,keepaspectratio]{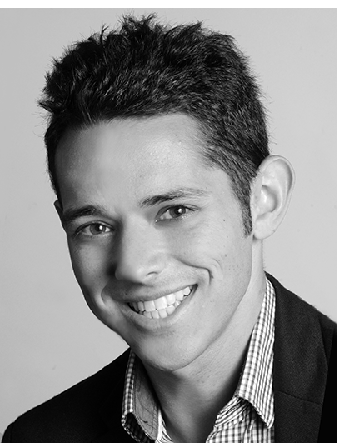}}]{Scott J. Moura}(S'09--M'13)
	received the B.S. degree from the University of California, Berkeley, CA, USA, and the M.S. and Ph.D. degrees from the University of Michigan, Ann Arbor, MI, USA, in 2006, 2008, and 2011, respectively, all in mechanical engineering.
	
	He is currently an Assistant Professor and Director of the Energy, Controls, and Applications Laboratory (eCAL) in Civil \& Environmental Engineering at the University of California, Berkeley. He is also an Assistant Professor  with the Smart Grid and Renewable Energy Laboratory, Tsinghua-Berkeley Shenzhen Institute. In 2011-2013, he was a postdoctoral fellow at the Cymer Center for Control Systems and Dynamics at the University of California, San Diego. In 2013 he was a visiting researcher in the Centre Automatique et Syst\`ems at MINES ParisTech in Paris, France. His current research interests include optimal and adaptive control, partial differential equation control, batteries, electric vehicles, and energy storage.
	
	Dr. Moura is a recipient of the National Science Foundation Graduate Research Fellowship, UC Presidential Postdoctoral Fellowship, O. Hugo Shuck Best Paper Award, ACC Best Student Paper Award (as advisor), ACC and ASME Dynamic Systems and Control Conference Best Student Paper Finalist (as student), Hellman Fellows Fund, University of Michigan Distinguished ProQuest Dissertation Honorable Mention, University of Michigan Rackham Merit Fellowship, College of Engineering Distinguished Leadership Award.
	
\end{IEEEbiography}

\end{document}